\documentclass[a4paper,12pt]{amsart}
\usepackage[cp1251]{inputenc}
\usepackage[T2A]{fontenc}
\usepackage[russian,ukrainian,english]{babel}
\usepackage{amscd}
\usepackage{amssymb}
\usepackage{amsmath,amsthm,amsfonts}
\usepackage{hhline}

\usepackage{cite}

\newtheorem{theorem}{Theorem}[section]
\newtheorem{lem}{Lemma}[section]

\newtheorem{ozn}{Definition}[section]
\newtheorem{remk}{Remark}[section]
\newtheorem{nas}{Corollary}[section]
\newtheorem{prop}{Proposition}[section]

\newcommand{\rr}{\mathbb{R}}

\renewcommand{\emptyset}{\varnothing}
\newcommand{\Cl}[1]{\overline{#1}}
\newcommand{\Int}[1]{\mathrm{Int}\,{#1}}

\newcommand{\cD}{\mathcal{D}}

\newcommand{\ff}{\mathfrak{f}}

\newcommand{\restrict}[2]{#1\raisebox{-0.65ex}{$\left|\vphantom{{#1}_{#2}}\right.$}_{#2}}

\usepackage{graphicx}

\begin{document}

\begin{center}
{\large\bf On the realization of graph as invariant of pseudoharmonic functions}
\end{center}
\vspace{3mm}
\begin{center}
{\small Polulyakh E.,\;Yurchuk I.}
\end{center}
\begin{center}
{\small Institute of Mathematics of Ukrainian National Academy of Sciences, Kyiv}
\end{center}
\vspace{5mm}

{\noindent \small\textbf{Abstract.} Necessary and sufficient conditions for a finite connected graph with a strict partial order on vertices to be a combinatorial invariant of pseudoharmonic function are obtained.}
\medskip


\section{Introduction}

This paper is a final step of research of pseudoharmonic functions~\cite{Kp,Yu,Bo,Mrt} defined on $D^2$ (i.e., two-dimensional closed disk). In~\cite{PY1} the invariant of such functions was constructed. It is a finite connected graph with a strict partial order and partial orientation on its vertices generated by a function. This invariant has several properties that follow from the nature of pseudoharmonic functions. In terms of such invariant the necessary and sufficient condition for functions to be topologically equivalent was obtained.

It says that two pseudoharmonic functions $f$ and $g$  are topologically equivalent iff there exists an isomorphism of combinatorial diagrams $\varphi:P(f)\rightarrow P(g)$ which preserves a strict partial order defined on them and the orientation (Th.3.1~\cite{PY1}).

In this paper, for simplicity, we will disregard the orientation and consider the combinatorial diagram \textit{\textbf{without orientation}}.

We will research the conditions for a finite connected graph $G$ with a strict partial order on its vertices to be a combinatorial invariant of some pseudoharmonic function.

In Section 2 the main theorem of the realization will be formulated. Let us comment Conditions A1-A3 that are necessary for Theorem~\ref{main_theorem}. $Cr$ -- cycle corresponds on a combinatorial diagram of a function $f$ to the restriction of $f$ onto the boundary of $D^2$. In~\cite{Kp} author proved that the connected components of level curves of pseudoharmonic function are isomorphic to a disjoin union of trees, so, Condition A2 is an analog of this statement. Finally, Condition A3 makes some vertices of $Cr$ -- cycle to be local extrema of function.

\begin{figure}[htbp]
\begin{center}
\includegraphics{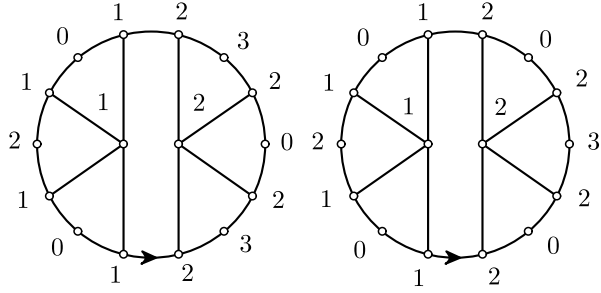}
\caption{On the left the graph can be realizated as diagram but on the right can not.}\label{fig_example}
\end{center}
\end{figure}

Let us remind some definitions from graph theory and previous results connected to these subjects.

Let $T$ be a tree (i.e., a $CW-$complex such that it does not contain cycles). Denote by $V$ the set of its vertices and by $V_{ter}(T)$ the subset of $V$ consisting of all terminal vertices. A disjoint union of trees is called a forest and denoted by $F=\bigcup\limits_i T_i$, where every $T_i$ is a tree. By a path $P(v_{i},v_{k})$  which connects vertices $v_{i}$ and $v_{k}$ of a graph $G$ we mean a sequence of adjacent edges $e_{j}$ such that each of them belongs to it once. It is known that for any two vertices of a tree there is a unique path connecting them. Two vertices of a graph are \textit{adjacent} if they are the ends of the same edge. A graph $G$ is embedded into $R^{2}$ if there exists a bijection $\phi:G\rightarrow \mathbb{R}^2$ such that two points $\phi(x)$ and $\phi(y)$ are joined by a segment iff $x$ and $y$ are joined by an edge of $G$ and no two distinct open segments have a point in common.

Let us remind necessary definitions and facts from~\cite{PY}.
Let $V^{\ast}$ be a subset of $V$ such that $V_{ter} \subseteq V^{\ast}.$

Let also $\varphi : T \rightarrow \rr^{2}$ be an embedding such that
\begin{equation}\label{eq_02}
\varphi(T) \subseteq D^{2} \,, \quad \varphi(T) \cap \partial D^{2} = \varphi(V^{\ast}) \,.
\end{equation}

\begin{lem}[see \cite{PY}]\label{lemma_01}
A set $\rr^{2} \setminus (\varphi(T) \cup \partial D^{2})$ has a finite number of connected components
\[
U_{0} = \rr^{2} \setminus D^{2}, U_{1}, \ldots, U_{m} \,,
\]
and for every $i \in \{1, \ldots, m\}$ a set $U_{i}$ is an open disk and is bounded by a simple closed curve
\[
\partial U_{i} = L_{i} \cup \varphi(P(v_{i}, v_{i}')) \,, \quad
L_{i} \cap \varphi(P(v_{i}, v_{i}')) = \{ \varphi(v_{i}), \varphi(v_{i}') \}
\]
 where $L_{i}$ is an arc of $\partial D^{2}$ such that the vertices $\varphi(v_{i})$ and $\varphi(v_{i}')$ are its endpoints, and  $\varphi(P(v_{i}, v_{i}'))$ is an image of the unique path $P(v_{i}, v_{i}')$ in $T$ which connects $v_{i}$ and $v_{i}'$.
\end{lem}

\begin{nas}\label{nas_01}
Let $T$ be a tree with a fixed subset of vertices $V^{\ast} \supseteq V_{ter}$ and $\varphi : T \rightarrow \rr^{2}$ an embedding which satisfies~\eqref{eq_02}.

Then the following conditions hold true.

1) In notation of Lemma~\ref{lemma_01}
\[
L_{i} \cap \varphi(T) = \{\varphi(v_{i}), \varphi(v_{i}')\} \,,
\quad i = 1, \ldots, m \,.
\]

2) If there exists an arc $L$ of circle $\partial D^{2}$ with the ends $\varphi(u_{1})$, $\varphi(u_{2})$ such that $L \cap \varphi(T) = \{\varphi(u_{1}), \varphi(u_{2})\}$ for some $u_{1}$, $u_{2} \in V^{\ast}$, then  there exists $k \in \{1, \ldots, m\}$ such that $L \cup \varphi(P(u_{1}, u_{2})) = \partial U_{k}$ (then $L = L_{k}$, $u_{1} = v_{k}$, $u_{2} = v_{k}'$).
\end{nas}

We also should remind that a ternary relation $O$ on a set $A$ is any subset of the $3^{rd}$ cartesian power $A^3:O\subseteq A^3$.

Let $A$ be a set, $O$ be a ternary relation on $A$ which is asymmetric ($(x,y,z)\in O \Rightarrow (z,y,x)\overline{\in} O $), transitive $(x,y,z)\in O, (x,z,u)\in O \Rightarrow (x,y,u) \in O$ and cyclic $(x,y,z)\in O \Rightarrow (y,z,x)\in O$. Then $O$ is called a cyclic order on the set $A$~\cite{Nov}.

A cyclic order $O$ is complete on a finite set $A$, $\sharp A \geq 3$, if for every $x,y,z\in A, x\neq y \neq z \neq x $ there exists a permutation $u,v,w$ of sequence $(x,y, z)$ such that $(u,v,w) \in O$.

\begin{prop}[see \cite{PY}]
Let there is a complete cyclic order $O$ on some finite set $A$, $\sharp A \geq 3$.

Then for every $a \in A$ there exist unique $a'$, $a'' \in A$ such that
\begin{itemize}
    \item $O(a', a, b)$ for all $b \in A \setminus \{a, a'\}$;
    \item $O(a, a'', b)$ for all $b \in A \setminus \{a, a''\}$,
\end{itemize}
and $a' \neq a''$.
\end{prop}

\begin{ozn}[see \cite{PY}]
Let there is a complete cyclic order $O$ on a set $A$, $\sharp A \geq 3$. Elements $a_{1}$, $a_{2} \in A$ are said to be \emph{adjacent} with respect to a cyclic order $O$ if one of the following conditions holds:
\begin{itemize}
    \item $O(a_{1}, a_{2}, b)$ for all $b \in A \setminus \{a_{1}, a_{2}\}$;
    \item $O(a_{2}, a_{1}, b)$ for all $b \in A \setminus \{a_{1}, a_{2}\}$.
\end{itemize}
\end{ozn}

\begin{remk}
From the previous Proposition it follows that every element has exactly two adjacent elements on a finite set $A$ with a complete cyclic order.
\end{remk}

We assume that if $\sharp V^{\ast} \geq 3$ then there is some cyclic order $C$ defined on $V^{\ast}$.

\begin{ozn}[see \cite{PY}]\label{ozn_01}
 A tree $T$ is called \emph{$\cD$-planar} if there exists an embedding $\varphi : T \rightarrow \rr^{2}$ which satisfies~\eqref{eq_02} and if $\sharp V^{\ast} \geq 3$ then a cyclic order $\varphi(C)$ on $\varphi(V^{\ast})$ coincides with a cyclic order which is generated by the orientation of $\partial D^{2} \cong S^{1}$.
\end{ozn}
The following theorem is the criterion of $\cD$-planarity of a tree.

\begin{theorem}[see \cite{PY}]\label{crit_d_planarity}
If $V^{\ast}$ contains just two vertices, a tree $T$ is $\cD$-planar.

If $\sharp V^{\ast} \geq 3$ then a $\cD$-planarity of $T$ is equivalent to satisfying the following condition:
\begin{itemize}
    \item for any edge $e$ there are exactly two paths such that they pass through an edge $e$ and connect two vertices of $V^{\ast}$ adjacent with respect to a cyclic order $C$.
\end{itemize}
\end{theorem}

Let us remind that a relation $<$ is called to be strict partial order if it is transitive, antireflexive and antisymmetric~\cite{Mel}.

\begin{ozn}[see \cite{PY}]\label{ozn_zruchne}
Let $A$ be a finite set. A binary relation $\rho$ on $A$ is said to be \emph{convenient} if
\begin{itemize}
    \item[1)] for all $a$, $b \in A$ from $a \rho b$ it follows that $a \neq b$;
    \item[2)] for every $a \in A$ there is no more than one $a' \in A$ such that $a \rho a'$;
    \item[3)] for every $a \in A$ there is no more than one $a'' \in A$ such that $a'' \rho a$.
\end{itemize}
\end{ozn}
Let $\mu$ be some relation on a set $A$.

\begin{ozn}[see \cite{PY}]\label{ozn_cycle}
Elements $b_{0}, \ldots, b_{n} \in A$, $n \geq 1$ are said to generate \emph{$\mu$-cycle} if a graph of the relation $\mu$ contains a set
\begin{equation}\label{eq_cycle}
\{(b_{0}, b_{1}), \ldots, (b_{n-1}, b_{n}), (b_{n}, b_{0})\} \,.
\end{equation}
\end{ozn}

\begin{ozn}[see \cite{PY}]\label{ozn_chain}
Elements $b_{0}, \ldots, b_{n} \in A$, $n \geq 0$ generate \emph{$\mu$-chain} if for arbitrary $a \in A$ the pairs $(a, b_{0})$ and $(b_{n}, a)$ do not belong to a graph of $\mu$ and for $n \geq 1$ a graph of the relation $\mu$ contains a set
\begin{equation}\label{eq_chain}
\{(b_{0}, b_{1}), \ldots, (b_{n-1}, b_{n})\} \,.
\end{equation}
\end{ozn}

\section{Conditions for the realization. Main theorem.}

Let $G \subset R^{3}$ be a finite connected graph with a strict partial order on vertices. We assume that every vertex of $G$ has a degree not less than 2.

A set $V\times V$ is divided onto two classes $C_{1}$ and $C_{2}$. Vertices $v_{1}$ and $v_{2}$ are contained in $C_{1}$ if they are comparable (i.e. either $v_{1}<v_{2}$ or $v_{2}<v_{1}$ holds true) and $C_{2}$ otherwise.

\begin{ozn}\label{ozn_cr_cycle}
\emph{Cr -- cycle} of $G$ is a subgraph $\gamma$ which is a simple cycle such that every pair of adjacent vertices of $\gamma$ belongs to $C_{1}$.
\end{ozn}

In what follows we will consider the following conditions on a graph $G\subset R^{3}$:
\begin{enumerate}
\item[A1)] there exists the unique $Cr$ -- cycle $\gamma$;

\item[A2)]$\overline{G\setminus \gamma}=F=\bigcup\limits_{i=1}^{k} T_{i}$, where $F$ is a forest such that
\begin{itemize}
\item if $v_{k}<v$ ($v_{k}>v$) for some vertex $v_{k}\in T_{i}\subset F$, where $v\in G$, then $v_{l}<v$ ($v_{l}>v$) for an arbitrary $v_{l}\in T_{i}\subset F$, $l\neq k$;

\item $deg(v)=2s\geq 4$ for an  arbitrary vertex $v\in G\setminus \gamma$;

\end{itemize}

\item[A3)] \textit{The condition for a strict order on  $Cr$ -- cycle $\gamma$}:
for any vertex $v$ of the subgraph $\gamma$ and its adjacent vertices $v_{1}$ and $v_{2}$ such that $v_{1}, v_{2}\in \gamma$ the following conditions hold true:
\begin{itemize}
 \item if $deg(v)=2$, then $deg(v_{1})>2$, $deg(v_{2})>2$ and there exists the unique index $i$ such that $v_{1},v_{2}\in T_{i}$;
 \item if $deg(v)=2s>2$ ($deg(v)=2s+1$), then $v_{1}\lessgtr v \gtrless v_{2}$ ($v_{1}\lessgtr v \lessgtr v_{2}$).
\end{itemize}

\item[A4)] \textit{The condition for a strict order on $G$}: if $v',v''\in C_2$, then from $v>v'$ it follows that $v>v''$.
\end{enumerate}
We remark that from $A2$ it follows that all vertices of any connected component $T_i$ are pairwise non comparable.

If A2 holds true, then, obviously, there exists a nonempty subset of vertices $V^{\ast}$ of $F$ which contains a set $V_{ter}$ of all terminal vertices of $F$ such that  $V^{\ast}=V(F)\cap \gamma$. It is clear that the subset of vertices $V^{\ast}$ of $F$ is divided onto the subsets $V^{\ast}_{k}$
such that $V_{ter}(T_k)\subset V_{k}^{\ast}\subset T_k\subset F$ and $V^{\ast}=\bigcup\limits_i V_{i}^{\ast}$.

\begin{ozn}\label{def_plan}
A finite graph $G\subset R^{3}$ is called $\mathfrak{D}-$~planar if there exists a subgraph $\gamma$ and an embedding $\varphi : G \rightarrow D^{2}$ such that the following conditions hold true:
\begin{itemize}
	\item $\gamma$ is a simple cycle;
	\item $\Cl{G \setminus \gamma} = \bigcup_{i=1}^{k} T_{i} = F$  is a finite union of trees;
	\item $\gamma$ contains all terminal vertices of $F$;
	\item $\varphi(\gamma) = \partial D^{2}$, $\varphi(G \setminus \gamma) \subseteq \Int{D^{2}}$.
\end{itemize}

\end{ozn}

\begin{theorem}\label{forest}
Let $G$ is a graph, $\gamma\subseteq G$ is a cycle such that $\Cl{G\setminus \gamma}=\bigsqcup_i T_i$, where every $T_i$ is a tree.

Then $G$ is $\mathfrak{D}$ -- planar if and only if every tree $T_{i}$ with the subset of vertices $V_i^\ast$ which has a cyclic order induced from $\gamma$ is $\mathfrak{D}$-planar and for any indexes $m$ and $n$ the subset of vertices $V_{n}^{\ast}$ of the tree $T_{n}$ belongs to a unique connected component of a set $\gamma\setminus V^{\ast}_{m}$, where $m\neq n$, $V^{\ast}_{j}\subset T_{j}$, $j=m,n$.
\end{theorem}

\begin{proof}
{\itshape\bfseries Necessity.} Suppose that a graph $G$ is $\mathfrak{D}$ -- planar. It is clear that every tree $T_{i}$ is $\mathfrak{D}$-planar. Let us assume that there exist some indexes $s$ and $l$ such that the subset of vertices $V_{s}^{\ast}$ of the tree $T_{s}$ belongs to two connected components $S'$ and $S''$ of the set $\gamma\setminus V^{\ast}_{l}$, where $V^{\ast}_{l}\subset T_{l}$. Let us assume that the subset $V_{s_{1}}^{\ast}$ of $V_{s}^{\ast}$ belongs to $S'$ and $V_{s_{2}}^{\ast}$ belongs to $S''$. Let us consider the following paths: $P_{1}$ which connects the ends of arc $S'$ ($S''$) (by construction they belong to $V^{\ast}_l$) in $G$ and $P_{2}$ which connects arbitrary two vertices of $V^{\ast}_{s}$ of the tree $T_{s}$ such that one of them belongs to $S'$, another belongs to $S''$. By our initial assumption there exists an embedding $\varphi$ of paths $\varphi(P_{1})$ and $\varphi(P_{2})$. They can be considered as two hordes which are contained into $Int D^2$ with the ends on $\partial D^2$. They have a common point which is not a vertex of $F$ ($T_{k}$ and $T_{l}$ are disconnected) since one pair of the ends parts another. It contradicts to the fact that $\varphi$ is an embedding.

\medskip

{\itshape\bfseries Sufficiency} will be proved by an induction on the number $n$ of trees in the forest $F = \Cl{G \setminus \gamma} = \bigsqcup_{i=1}^n T_i$.

Let us regard $G$ as a cell complex. Then the cycle $\gamma$ considered as a subspace of the topological space $G$ is homeomorphic to a circle. Fix an orientation on $\gamma$. It induces a cyclic order on it. Now we induce from $\gamma$ a cyclic order on each $V^{\ast}_i$ with $\sharp V^{\ast}_i > 2$, $i = 1, \ldots, n$.

We should remark that the following is straightforward: if for a fixed orientaition of $\gamma$ a tree $T_i$ with the cyclic order on $V^{\ast}_i$ induced from $\gamma$ is $\mathfrak{D}$-planar, then for an inverse orientation of $\gamma$ a tree $T_i$ with the cyclic order on $V^{\ast}_i$ induced from that orientation of $\gamma$ is also $\mathfrak{D}$-planar. So the choice of an orientation of $\gamma$ does not affect the ongoing considerations.

Suppose that every tree $T_{k}$ is $\mathfrak{D}$-planar and for arbitrary indexes $r$ and $s$ the subset $V_{s}^{\ast}$ of vertices of $T_{s}$ belongs to the unique connected component of the set $\gamma\setminus V^{\ast}_{r}$, where $r \neq s$, $V^{\ast}_{j} = T_{j} \cap \gamma$, $j=r,s$.

\emph{Basis of induction}. Let $F = T_1$.

First let $V^{\ast}_1 = \{v_1, v_2\}$ for some $v_1, v_2 \in V(T_1)$. Since $T_1$ is $\mathfrak{D}$-planar, then there is an embedding $\varphi_1 : T_1 \rightarrow D^2$ such that $\varphi(T_1) \cap \partial D^2 = \varphi(V^{\ast}_1) = \{\varphi(v_1), \varphi(v_2)\}$. Obviously, the cycle $\gamma$ consists of two edges with common endpoints $v_1$ and $v_2$. Fix some embedding $\varphi' : \gamma \rightarrow \partial D^2$ such that $\varphi_1(v_i) = \varphi'(v_i)$, $i = 1, 2$. Now it is straightforward to see that the mapping $\varphi : G \rightarrow D^2$,
\begin{equation}\label{eq_embedding_1}
\varphi(\tau) =
\begin{cases}
\varphi_1(\tau), & \mbox{if $\tau$ is in $T_1$,}\\
\varphi'(\tau), & \mbox{when $\tau$ is in $\gamma$,}
\end{cases}
\end{equation}
is well defined and complies with Definition~\ref{def_plan}.

Now let $\sharp V^{\ast}_1 > 2$.
At first we define a bijective and continuous map $\varphi':\gamma\rightarrow \partial D^{2}$ such that the cyclic order induced on $\varphi'(\gamma)$ by $\varphi'$ coincides with the cyclic order induced by the positive orientation of $\partial D^2$. It is obvious that it is an embedding and  $\varphi'(\gamma)=\partial D^{2}$.

Let us consider the tree $T_{1}$. From its $\mathfrak{D}$-planarity it follows that there exists an embedding $\varphi_{1} : T_{1} \rightarrow D^{2}$ such that $\varphi_{1}(T_{1}) \subset D^{2}$, $\varphi_{1}(V^{\ast}_{1}) \subset \partial D^{2}$,
$\varphi_{1}(T_{1} \setminus V^{\ast}_{1}) \subset Int D^{2}$, where $V_{ter}(T_{1})\subseteq V^{\ast}_{1} \subset V$.
We can choose $\varphi_1$ in such way that $\varphi_{1}|_{V^{\ast}_{1}} = \varphi'|_{V^{\ast}_{1}}$ since a cyclic order on vertices of $V^{\ast}_1$ is consistent with the cyclic order on $\partial D^{2} = \varphi'(\gamma)$ which is in turn induced by $\varphi'$ from the cyclical order on $\gamma$.

Then it is easy to see that the mapping $\varphi : G \rightarrow D^2$ given by~\eqref{eq_embedding_1} is well defined and satisfies all requirements of Definition~\ref{def_plan}.

\emph{Step of induction.} Let $G = \gamma \cup F$, $F = \Cl{G \setminus \gamma} = \bigsqcup_{j=1}^n T_j$, $n > 1$. Suppose that for any graph $G'$ with a cycle $\gamma'$ such that $F' = \Cl{G' \setminus \gamma'} = \bigsqcup_{i=1}^k T_i'$ is a forest and $k < n$ our Theorem holds true.

First we are going to prove that there is a tree $T_s$ in $F$ such that the set $\bigcup_{j \neq s} V^{\ast}_{j}$ is contained in a single connected component of the set $\gamma\setminus V^{\ast}_{s}$.

For every $i = 1, \ldots, n$ we shall denote by $\nu(T_i)$ the maximal cardinality of subsets $M_i \subseteq \{1, \ldots, n\}$ which satisfy the following property: a set $\bigcup_{j \in M_i} V^{\ast}_{j}$ is contained in a single connected component of the set $\gamma\setminus V^{\ast}_{i}$.

As $n > 1$ then $1 \leq \nu(T_i) \leq n-1$, $i = 1, \ldots, n$. And $\nu(T_i) = n-1$ iff a set $\bigcup_{j \neq i} V^{\ast}_{j}$ is contained in a single connected component of the set $\gamma\setminus V^{\ast}_{i}$.

Let $\nu(T_i) < n-1$ for a certain $i$. Let us designate by $\gamma^i_1, \ldots, \gamma^i_{m(i)}$ all components of the complement $\gamma \setminus V^{\ast}_i$. Then there exist at least two different indexes $r'$ and $r''$ for which relations $\gamma^i_{r'} \cap \bigcup_{j \neq i} V^{\ast}_j \neq \emptyset$ and $\gamma^i_{r''} \cap \bigcup_{j \neq i} V^{\ast}_j \neq \emptyset$ hold true.

We can select $r'$ in such way that $\gamma^i_{r'} \cap \bigcup_{j=1}^n V^{\ast}_j = \gamma^i_{r'} \cap \bigcup_{j \in M_i} V^{\ast}_j$ for a subset $M_i$ of $\{1, \ldots, n\} \setminus \{i\}$ with cardinality $\sharp M_i = \nu(T_i)$. Fix $i' \notin M_i \cup \{i\}$ and let $r''$ be an index such that $V^{\ast}_{i'} \subset \gamma^i_{r''}$. It is clear that $r' \neq r''$.
Since both $V^{\ast}_i$ and $\bigcup_{j \in M_i} V^{\ast}_j$ are contained in a connected subset $\gamma \setminus \gamma^i_{r''}$ of the cycle $\gamma$ and $(\gamma \setminus \gamma^i_{r''}) \cap V^{\ast}_{i'} \subset (\gamma \setminus \gamma^i_{r''}) \cap \gamma^i_{r''} = \emptyset$ then the set $V^{\ast}_i \cup \bigcup_{j \in M_i} V^{\ast}_j$ lies in a single component of the complement $\gamma \setminus V^{\ast}_{i'}$ and consequently $\nu(T_{i'}) \geq \nu(T_i) + 1$.

So, in a finite number of steps we shall find an index $s$ such that $\nu(T_s) \geq n-1$, therefore a set $\bigcup_{j \neq s} V^{\ast}_{j}$ is contained in a single connected component of the set $\gamma\setminus V^{\ast}_{s}$.

Without loss of generality we can regard that $\nu(T_n) = n-1$. Repeating the argument we used to verify the base of induction we can find an embedding $\varphi_n : \gamma \cup T_n \rightarrow D^2$ which maps $\gamma$ onto $\partial D^2$ and such that an orientation on $\varphi_n(\gamma) = \partial D^2$ induced by $\varphi_n$ coincides with the positive orientation on this set induced from $D^2$.

Lemma~\ref{lemma_01} implies that $D^2 \setminus \varphi_{n}(T_{n})= \bigcup_s U_s$, where $\overline{U}_{s}\cong D^2$ and $\partial \overline{U}_{s} \subset \partial D^{2} \cup \varphi_{n}(T_{n})$ for any $s$. By the choice of $T_n$ the subset $\bigcup_{i \neq n} V^{\ast}_{i}$ of vertices of a forest $F' = \bigcup_{i = 1}^{n-1} T_i$ belongs to a single connected component $\gamma_0$ of the set $\gamma \setminus V^{\ast}_n$. From this and from Corollary~\ref{nas_01} it follows that there exists an index $m$ such that a domain $U_{m}$ satisfies the inclusions $\varphi_n(\gamma_0) \subset (\Cl{U}_m \cap \partial D^2)$, $\partial U_m \setminus \varphi_n(\gamma_0) = \varphi_n(P)$, where $P = P(v', v'')$ is a path in $T_n$ which connects a pair of vertices $v'$, $v'' \in V^{\ast}_n$.

Let us consider a cycle $\gamma' = \gamma_0 \cup P$ in $G$. It is clear that it is simple. Denote
\[
G' = \gamma' \cup F' = \gamma' \cup \bigcup_{i=1}^{n-1} T_i \,.
\]
Since $F' \cap \gamma = \bigcup_{i=1}^{n-1} V^{\ast}_i \subset \gamma_0$ by construction, then $F' \cap \gamma' \subset \gamma_0$ and $F' = \Cl{G' \setminus \gamma'}$.

The following claim is straightforward.
\emph{Suppose we have two oriented cicrles $S_1$ and $S_2$ and two arcs $\gamma_1 \subset S_1$ and $\gamma_2 \subset S_2$ such that orientaion of each arc is coordinated with an orientation of the corresponding circle. Let $\Phi : \gamma_1 \rightarrow \gamma_2$ be an orientation preserving homeomorphism. Let also $O_k$, $k=1,2$, be a full cyclic order on $S_k$ induced by the orientation of $S_k$. Then $O_2|_{\gamma_2} = \Phi(O_1|_{\gamma_1})$.}

Let us induce an orientation on $\gamma_0$ from $\gamma$ and choose an orientation on $\gamma'$ which is coordinated with the selected orientation of $\gamma_0$. Let $\Phi = Id : \gamma_0 \rightarrow \gamma_0$. Then by the claim above cyclic orders on $\gamma_0$ induced from $\gamma$ and from $\gamma'$ should coincide.

Every tree $T_i$, $i \in \{1, \ldots, n-1\}$, with the subset of vertices $V^{\ast}_i$ which has a cyclic order induced from the positive orientation of $\gamma$ is $\mathfrak{D}$-planar by our initial assumption. As $V^{\ast}_i \subset \gamma_0$, then according to what was said above every $T_i$ is $\mathfrak{D}$-planar with respect to a cyclic order induced on $V^{\ast}_i$ from positive orientation of $\gamma'$.

It is easy to see that since the set $\bigcup_{i=1}^{n-1} V^{\ast}_i$ is contained in the connected set $\gamma_0 \subseteq \gamma \cap \gamma'$ and by our initial assumption for any indexes $j, k \in \{1, \ldots, n-1\}$ the subset of vertices $V_{j}^{\ast}$ of the tree $T_{j}$ belongs to a unique connected component of a set $\gamma \setminus V^{\ast}_{k}$, where $j \neq k$, then every set $V_{j}^{\ast}$ is contianed in a single connected component of a set $\gamma' \setminus V^{\ast}_{k}$, $j \neq k$, $j, k \in \{1, \ldots, n-1\}$.

As a consequence from said graph $G' = \gamma' \cup \bigcup_{i=1}^{n-1} T_i$ is $\mathfrak{D}$-planar by the inductive hypothesis. So, there exists an embedding $\varphi' : G' \rightarrow D^2$ which is compliant with Definition~\ref{def_plan}.

Then $\varphi'(\gamma') = \partial D^2$. Let us remind that by construction we have $\varphi_n(\gamma') = \partial U_m$. Evidently, a map $\psi_0 = \varphi_n \circ (\varphi')^{-1} : \partial D^2 \rightarrow \partial U_m$ is homeomorphism.
Let us remind (see \cite{Newman}) that every homeomorphism of simple closed curves in the plane can be extended to a homeomorphism of disks bounded by these curves. So, there exists a homeomorphism $\psi : D^2 \rightarrow \Cl{U}_m$ such that $\psi|_{\partial D^2} = \psi_0$.

Let us consider a map $\varphi : G \rightarrow D^2$ defined by the relation
\[
\varphi(\tau) =
\begin{cases}
\varphi_n(\tau)\,, & \mbox{when $\tau \in \gamma \cup T_n$}\,,\\
\psi \circ \varphi'(\tau)\,, & \mbox{if $\tau \in T_i$, $i \in \{1, \ldots, n-1\}$}\,.
\end{cases}
\]
Since $(\gamma \cup T_n) \cap (\bigcup_{i=1}^{n-1} T_i) \subseteq \gamma'$ by construction and $\psi \circ \varphi'(\tau) = \psi_0 \circ \varphi'(\tau) = \varphi_n(\tau)$ for every $\tau \in \varphi'(\gamma') = \partial D^2$, then $\varphi$ is well defined. The sets $\gamma \cup T_n$ and $\bigcup_{i=1}^{n-1} T_i$ are closed, so $\varphi$ is continuous. And it is straightforward to see that this map is injective. Therefore $\varphi$ is the embedding of compact $G$ into $D^2$.

By our initial assumptions every tree $T_k$, $k=1, \ldots, n$, is $\mathfrak{D}$-planar with respect to the cyclic order on the set $V^{\ast}_k = T_k \cap \gamma$ induced from $\gamma$. Then $V_{ter}(T_k) \subseteq V^{\ast}_k \subset \gamma$, hence $\gamma$ contains all terminal vertices of the forest $F$.

Finally, observe that $\varphi(\gamma) = \varphi_n(\gamma) = \partial D^2$.

So, graph $G$ satisfies all conditions of Definition~\ref{def_plan} and by induction principle conditions on $G$ to be $\mathfrak{D}$-planar stated in Theorem are sufficient for $G$ with any number of trees in a forest $F = \Cl{G \setminus \gamma}$.
\end{proof}

\begin{remk}\label{for_nas}
If $G$ satisfies $A1$ and $A2$, then Theorem~\ref{forest} holds true for it.
\end{remk}

Assume that $G$ satisfies A1 and A2.
Let us consider arbitrary two vertices $v_{1},v_{2}$ of the set $V^{\ast}_{i}$ of the subgraph $T_{i}$ of $G$. The set $\gamma \setminus (v_{1}\cup v_{2})$ consists of disjoint union of two connected sets $\gamma_1$ and $\gamma_2$.

\begin{ozn} Pair of vertices $v_1,v_2\in V_i^\ast$ is called \emph{boundary} if either $\gamma_1$ or $\gamma_2$ does not contain any vertex of $V^{\ast}_{i}$ and at least one vertex of $V^{\ast}\setminus V^{\ast}_{i}$ belongs to it.
\end{ozn}
Denote by $\omega(v_{1},v_{2})$ the boundary pair, designate by $\alpha$ the set $\gamma_k$ which does not contain any vertex of $V^{\ast}_{i}$ and at least one vertex of $V^{\ast}\setminus V^{\ast}_{i}$ belongs to it. It is clear that for every vertex $v_j$ of the boundary pair $\omega(v_{1},v_{2})$ there exists an adjacent vertex $\tilde{v}_j$ such that $\tilde{v}_j \in \alpha$, where $j=\overline{1,2}$.

\begin{ozn} A graph $G$ is called \emph{special} if the following conditions hold true:
\begin{itemize}
 \item[S1)] $G$ satisfies A1 and A2;
 \item[S2)] $G$ is $\mathfrak{D}$ -- planar;
 \item[S3)] for arbitrary boundary pair $\omega(v_{1},v_{2})\in V_i^\ast$ the pair of adjacent vertices $\tilde{v}_1,\tilde{v}_2$ belongs to the unique set $V^{\ast}_{k}$, where $V^{\ast}_{k}\subset V^{\ast}\setminus V_i^{\ast}$, $\tilde{v}_1,\tilde{v}_2\in\alpha$.
\end{itemize}
\end{ozn}

\begin{remk}
If $v_1$, $v_2$ is a boundary pair, then the pair $\tilde{v}_1,\tilde{v}_2$ is the boundary pair of a tree $T_{k}\supset V^{\ast}_{k}$ for a special graph.
\end{remk}

\begin{lem}
\label{l2}
If a graph $G$ is special, then the set $\Theta=D^{2}\setminus \varphi(G)$ consists of disjoint union of the domains $U_{i}$ such that $\partial\overline{U}_{i}$ contains either one or two nondegenerate arcs of the boundary $\partial D^{2}$, where $\varphi:G\rightarrow D^{2}$ is an embedding such that $\varphi(\gamma)=\partial D^{2}$, $\varphi(G\setminus\gamma)\subset Int D^{2}$.
\end{lem}
\begin{proof}
Let $\varphi: \Gamma\rightarrow D^{2}$ be an embedding of the special graph $G$ such that $\varphi(\gamma)=\partial D^{2}$, $\varphi(G\setminus
\gamma)\subset Int D^{2}$. Condition A2 holds true hence there does not exist a domain $U_{i}$ such that its boundary $\partial
\overline{U}_{i}$ does not contain an arc of $\partial D^{2}$.
The set of vertices $\bigcup\limits_k\{\varphi(V^{\ast}_{k})\}$ divides $\partial D^{2}$ onto the arcs $A_{j}$. Let $\varphi(v_{i})$ and $\varphi(v_{i+1})$ be the end points of $A_{i}$. Let us consider two cases:

{\emph{Case 1:}} $v_{i},v_{i+1}\in T_{k}$. From Corollary~\ref{nas_01} it follows that $\partial\overline{U}_{i}$ contains one arc $\partial D^{2}$.

{\emph{Case 2:}} $v_{i}\in T_{n} ,v_{i+1}\in T_{m}.$
By moving along $\partial D^{2}\setminus A_i$ from $\varphi(v_{i+1})$ ($\varphi(v_{i})$) in the direction of $\varphi(v_{i})$ ($\varphi(v_{i+1})$) we find the first vertex $\varphi(v_{j})$ such that $v_{j}\in V^{\ast}_n\subset T_{n}$
($v_{j}\in V^{\ast}_m \subset T_{m}$). It is clear that there exists the unique path
$P(v_{i},v_{j})$ ($P(v_{i+1},v_{j})$) such that $P(v_{i},v_{j})\in
T_{n}$ ($P(v_{i+1},v_{j})\in T_{m}$). Condition $S2$ holds true hence the pair $v_{i},v_{j}$ $(v_{i+1},v_{j})$ is boundary and by $S3$ for $v_{j}$ there exists a vertex $v_{j-1}$ ($v_{j+1}$) adjacent to $v_j$
such that $v_{j-1}\in T_{m}$ ($v_{j+1}\in T_{n}$). Thus the domain $\overline{U}_{i}$ such that $\varphi(v_{i}),\varphi(v_{i+1}),\varphi(v_{j}),\varphi(v_{j-1})\in
\partial\overline{U}_{i}$ ($\varphi(v_{i}),\varphi(v_{i+1}),\varphi(v_{j}),\varphi(v_{j+1})\in
\partial\overline{U}_{i}$) contains two boundary arcs of $\partial D^{2}$.
\end{proof}

\begin{ozn}
A special graph $G\subset R^{3}$ is called $\Delta$ -- graph if it satisfies A3.
\end{ozn}

\begin{lem}
\label{l3}
If $\hat{v}=min \{V\}$, $\check{v}=max \{V\}$ are vertices of $\Delta$ -- graph $G$, then $\hat{v},\check{v}\in \gamma$ and $deg(\hat{v})=deg(\check{v})=2$.
\end{lem}

\begin{proof}
We prove lemma for the case of minimal value $\hat{v}=min \{V\}$. Without loss of generality, suppose that $\hat{v}\in T_{j}$, where $T_{j}$ is a tree. From A2 it follows that there exists some vertex $v'\in \gamma \bigcap V_{ter}(T_{j})$ ( $v'$ is terminal of $T_{j}$). Condition A3 holds true hence for $v'$ there exists an adjacent vertex $v_{1}$  such that $v_{1}< v'$. It contradicts to $A2$ since $\hat{v}=min V$. It follows that $ \hat{v}$ belongs to a set $\gamma\setminus \bigcup\overline{T_{i}}$ which contains only the vertices of degree 2.

The case $\check{v}=max \{V\}$ is proved similarly.
\end{proof}

Let us remind some definitions~\cite{RF}.

A cover $\Gamma$ of a space $X$ is called \emph{fundamental} if arbitrary set such that its intersection with any set $B\in \Gamma$ is open in $B$ is also open in $X$. All finite and locally finite closed covers are fundamental.

Let $\Gamma$ be a fundamental cover of $X$ and for any set $A\in \Gamma$ a continuous map $f_{A}:A\rightarrow Y$ is defined such that if $x\in A\cap B (A,B\in \Gamma)$ then $f_{A}(x)=f_{B}(x)$. It is known that a map $f:X\rightarrow Y$ defined as $f(x)=f_{A}(x)$, where $x\in A$, $A\in \Gamma$, is continuous.

Let $A$ be a finite set. It is obvious that a function $g : A \rightarrow \rr$ induces a partial ordering relation on the set $A$ by correlation
\[
a' < a'' \quad \mbox{if } g(a') < g(a'') \,, \quad a', a'' \in A \,.
\]

Suppose that there are two partial orders $<$ and $<'$ on $A$. We will say that \emph{a partial order $<'$ extends an order $<$} if the identical map $Id : (A, <) \rightarrow (A, <')$ is monotone.

\begin{lem}\label{main_lemma}
Let us consider $\Delta$-graph $G$ as $CW-$complex. There exists a continuous function $g : G \rightarrow \rr$ on the topological space $G$ which satisfies the following conditions:
\begin{itemize}
	\item $g$ maps a partially ordered set $V(G)$ of vertices of $G$ into $\rr$ monotonically;
	\item local extrema of the restriction $\restrict{g}{\gamma}$ are exactly vertices of $G$ with even degree which belong to the cycle $\gamma$;
	\item any tree $T_{i}$, $i = 1, \ldots, k$ of $F = \Cl{G \setminus \gamma}$ is contained in some level set of the function $g$.
\end{itemize}

Then a partial order $<'$ induced by $g$ on the set $V(G)$ of vertices of $G$ is an extension of a partial order $<$ on $V(G)$.
\end{lem}

\begin{proof}
Let us consider a partition $\ff$ of the set $V(G)$ elements of which are vertices with degree 2 (they belong to $\gamma \setminus F$ by $A2$) and sets $V(T_{i})$, $i = 1, \ldots, k$ of vertices of trees of $F$.

From Condition $A2$ easily follows that relation of partial order on the set $V(G)$ induces a partial order on the quotient set $\hat{V} = V(G) / \ff$. Let us denote a projection map by $\pi : V(G) \rightarrow \hat{V}$. It is monotone by the construction.

It is evident that there exists a monotone map $\hat{g} : \hat{V} \rightarrow \rr$. A composition $g = \hat{g} \circ \pi : V(G) \rightarrow \rr$ is a monotone map as a composition of monotone maps.
From the construction it follows that any set $V(T_{i})$, $i = 1, \ldots, k$ belongs to some level set of a function $g$.

For any edge $e \in E(G)$ we fix a homeomorphism $\hat{g}_{e} : e \rightarrow [0, 1]$. It evidently maps the endpoints of $e$ on the set $\{0, 1\}$, therefore $\hat{g}_{e}^{-1}(\{0, 1\}) \subseteq V(G)$, $g \circ \hat{g}_{e}^{-1}(0)$ and $g \circ \hat{g}_{e}^{-1}(1)$ are defined.

For every $t \in E(G)$ we consider a monotone function
\begin{align*}
h_{e} & : [0, 1] \rightarrow [\min (g \circ \hat{g}_{e}^{-1}(0), g \circ \hat{g}_{e}^{-1}(1)), \max (g \circ \hat{g}_{e}^{-1}(0), g \circ \hat{g}_{e}^{-1}(1))] \,, \\
h_{e} & : t \mapsto (1-t) g(\hat{g}_{e}^{-1}(0)) + t g(\hat{g}_{e}^{-1}(1)) \,,
\end{align*}
and also a map

\[
g_{e} = h_{e} \circ \hat{g}_{e} : e \rightarrow [\min (g \circ \hat{g}_{e}^{-1}(0), g \circ \hat{g}_{e}^{-1}(1)), \max (g \circ \hat{g}_{e}^{-1}(0), g \circ \hat{g}_{e}^{-1}(1))] \,.
\]

It is obvious that for any two edges $e_{1}$, $e_{2} \in E(G)$ which have a common endpoint $v \in V(G)$ it holds true that $g_{e_{1}}(v) = g_{e_{2}}(v) = g(v)$. This allows us to extend a function $g$ on the edges of $G$ with the help of the following correlation
\[
g(x) = g_{e}(x) \,, \quad \mbox{for } x \in e \,.
\]
The set of all edges of $G$ generates closed covering of a topological space $G$. Graph is finite therefore such covering is fundamental and $g : G \rightarrow \rr$ is continuous, see above.

It is also obvious that if $g(v') = g(v'')$ for endpoints $v'$, $v''$ of some edge $e$, then $g_{e}(e) = g(v') = g(v'')$. Otherwise, a map $g_{e}$ is a homeomorphism.
By the construction we get that $g(v') = g(v'')$, where $v'$, $v'' \in V(T_{i})$ and $T_{i} \in F$. Therefore $g(T_{i}) = c_{i} \in \rr$, $i = 1, \ldots, k$, and any tree $T_{i}$ belongs to some level set of $g$.

Next to the last condition of Lemma easily follows from $A3$.
\end{proof}

\begin{theorem}\label{main_theorem}
If a graph $G$ is a combinatorial diagram of some pseudoharmonic function $f$, then $G$ is $\Delta$ -- graph.

If a graph $G$ is $\Delta$ -- graph, then a partial order on $V(G)$ can be extended so that the graph $G$ with a new partial order on the set of vertices will be isomorphic to a combinatorial diagram of some pseudoharmonic function $f$.
\end{theorem}

\begin{proof}
In order to prove the first part we should show that for a diagram $P(f)$ of pseudoharmonic function $f$ Condition $S3$ holds true. Suppose that for some boundary pair $\omega(v_{1},v_{2})\in q(f)\cap T_i$ (the existence of which follows from $C1-C3$, see~\cite{PY1}) the adjacent pair of vertices $\tilde{v}_1,\tilde{v}_2\in q(f)$ belongs to different sets $T_k$ and $T_{l}$, where $i\neq k$,$i\neq l$, $l\neq k$ and $\tilde{v}_1\in T_{k}$, $\tilde{v}_2\in T_{l}$. Then for a vertex $\tilde{v}_1$ ($\tilde{v}_2$) there exists $\widetilde{\widetilde{v}}_1$ ($\widetilde{\widetilde{v}}_2$) such that $\widetilde{\widetilde{v}}_1\in q(f)\cap T_{k}$ ($\widetilde{\widetilde{v}}_2\in q(f)\cap T_{l}$) and the pair $\tilde{v}_1,\widetilde{\widetilde{v}}_1$ ($\tilde{v}_2,\widetilde{\widetilde{v}}_2$) is boundary for the tree $T_{k}$ ($T_{l}$). It means that a domain $U_i$ such that $\partial U_i \ni \varphi(v_{1}),\varphi(v_{2}),\varphi(\tilde{v}_2),\varphi(\tilde{v}_1),
\varphi(\widetilde{\widetilde{v}}_1),\varphi(\widetilde{\widetilde{v}}_2)$ contains more than two boundary arcs but it contradicts to Lemma~\ref{l2}.

\medskip
Let us prove the second part of theorem.
Suppose that a graph $G$ is $\Delta$ -- graph.
Then there exists an embedding $\varphi: G\rightarrow D^{2}$ such that $\varphi(\gamma)=\partial D^{2}$, $\varphi(G\setminus \gamma)\subset Int D^{2}$. From Lemma~\ref{l2} it follows that the set $\Theta=D^{2}\setminus \varphi(G)$ consists of disjoint union of domains $U_{i}$ such that $\partial\overline{U}_{i}$ contains either one or two arcs of boundary $\partial D^{2}$.

Next we fix a continuous function $g : G \rightarrow \rr$ that satisfies the conditions of Lemma~\ref{main_lemma} and consider a continuous function $f = g \circ \varphi^{-1} : \varphi(G) \rightarrow \rr$ on the set $\varphi(G)$.

Our aim is to extend $f$ on all $U_{i}$ in order to obtain a continuous function on $D^{2}$ which can be locally represented as a projection on coordinate axis in a neighborhood of every point of $\Theta$.

Let us consider two types of domains.

\medskip

\emph{Case 1:} Let $U_{k}\subset \Theta$ be a domain such that  $\partial\overline{U}_{k}$ contains only one boundary arc $\alpha \subset \partial D^{2}$ and $\partial U_{k}\setminus \alpha=\beta$, where $\beta \subset \varphi(F)$.

It is clear that the set $\beta$ is connected therefore there exists the tree $T_{i} \subseteq F$ such that $\beta \subseteq \varphi(T_{i})$. From Lemma~\ref{main_lemma} it follows that $\restrict{f}{\beta} = \restrict{g}{T_{i}} = const$. Let $f(\beta) = c_{i} \in \rr$.

Let us consider the arc $\alpha$ and the preimage $\varphi^{-1}(\Cl{\alpha}) \subseteq \gamma$ and denote by $y'$ and $y''$ the endpoints of an arc $\alpha$. Then $v' = \varphi^{-1}(y')$ and $v'' = \varphi^{-1}(y'')$ belong to the set $V(T_{i}) \cap V(\gamma)$. The vertices $v'$ and $v''$ can not be adjacent vertices of the cycle $\gamma$ since they belong the same tree $T_{i}$ of $F$. Therefore $v'$ and $v''$ are non comparable, see Definition~\ref{ozn_cr_cycle} and Condition $A2$. Thus the set $\varphi^{-1}(\alpha)$ contains at least one vertex of the graph $G$ except $v'$ and $v''$. It is obvious that the arc $\alpha$ can not contain images of vertices of $F$ besides its endpoints $y'$ and $y''$. Therefore from $A3$ it follows that the arc $\alpha \setminus \{y', y''\}$ contains an image of exactly one vertex $y = \varphi(v)$ of degree 2. From Lemma~\ref{main_lemma} it follows that $f$ has a local extremum in the point $y$ and the arc $\alpha \setminus \{y', y''\}$ does not contain another local extrema of $f$. Suppose that $f(y) = c$. Points $v$ and $v'$ are comparable since $c \neq c_{i}$.

Denote by $\alpha'$ and $\alpha''$ subarcs of the arc $\alpha$. Suppose that the first of them connects points $y'$ and $y$ and the second connects $y$ and $y''$. From the above discussion it follows that $f$ is monotone on both arc $\alpha'$ and $\alpha''$. Thus maps $\psi' : \alpha' \rightarrow [0, 1]$, $\psi'' : \alpha'' \rightarrow [0, 1]$,
\begin{align*}
\psi'(z) & = \frac{f(z)-f(y)}{f(y')-f(y)} = \frac{f(z)-c}{c_{i}-c} \,, \\
\psi''(z) & = \frac{f(z)-f(y)}{f(y'')-f(y)} = \frac{f(z)-c}{c_{i}-c} \,,
\end{align*}
are homeomorphisms, in addition, $\psi'(y) = \psi''(y) = 0$, $\psi'(y') = \psi''(y'') = 1$.

\begin{figure}[htbp]
\begin{center}
\includegraphics{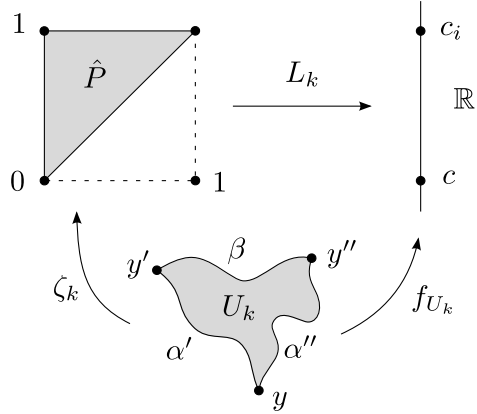}
\caption{Function on a simple connected domain with one boundary arc.}\label{fig_case_1}
\end{center}
\end{figure}

Let us consider a set $\hat{P} = \{ (x, y) \in [0, 1]^{2} \;|\; y \geq x \}$ and a map $L_{k} : \hat{P} \rightarrow \rr$,
\[
L_{k} : (x, y) \mapsto c (1-y) + c_{i} y \,.
\]
It is obvious that $L_{k}(0, 0) = c = f(y)$, $L_{k}([0, 1] \times \{1\}) = c_{i}$.

The arc $\beta$ is obviously homeomorphic to segment. Let $\psi : \beta \rightarrow [0, 1]$ be a homeomorphism such that $\psi(y') = 0$, $\psi(y'') = 1$.
We consider a map $\zeta_{k}^{0} : \partial U_{k} \rightarrow \partial \hat{P}$,
\[
\zeta_{k}^{0}(z) =
\left\{
\begin{array}{ll}
	(0, \psi'(z))\,, & \mbox{for } z \in \alpha' \,, \\
	(\psi''(z), \psi''(z))\,, & \mbox{for } z \in \alpha'' \,, \\
	(\psi(z), 1)\,, & \mbox{for } z \in \beta \,.
\end{array}
\right.
\]

It is easy to show that $\zeta_{k}^{0}$ is homeomorphism. Both the sets $\partial U_{k}$ and $\partial \hat{P}$ are simple closed curves thus we can use Schoenflies's theorem~\cite{Newman} and extend the homeomorphism $\zeta_{k}^{0}$ to $\zeta_{k} : \Cl{U}_{k} \rightarrow \hat{P}$.

Let us consider a continuous function, see Fig.~\ref{fig_case_1}
\[
f_{U_{k}} = L_{k} \circ \zeta_{k} : \Cl{U}_{k} \rightarrow \rr \,.
\]
It is obvious that this function locally  can be represented as a projection on coordinate axis in all points of $\Cl{U}_{k} \setminus \{y, y', y''\}$ .

let us prove that $\restrict{f_{U_{k}}}{\partial U_{k}} = \restrict{f}{\partial U_{k}}$. Really, for any $z \in \beta$ we have
\[
f_{U_{k}}(z) = L_{k} \circ \zeta_{k}(z) = L_{k}(\psi(z), 1) = c_{i} = f(z) \,;
\]
for $z \in \alpha'$ the following relations hold true
\begin{align*}
f_{U_{k}}(z) & = L_{k}(0, \psi'(z)) = c (1-\psi'(z)) + c_{i} \psi'(z) = \\
	& = \frac{c_{i}-f(z)}{c_{i}-c} \cdot c + \frac{f(z)-c}{c_{i}-c} \cdot c_{i} = f(z) \,;
\end{align*}
similarly, for $z \in \alpha'$ we have
\[
f_{U_{k}}(z) = L_{k}(\psi''(z), \psi''(z)) = c (1-\psi''(z)) + c_{i} \psi''(z) = f(z) \,.
\]

\medskip

\emph{Case 2:} Let $U_{k}\subset \Theta$ be a domain such that $\partial \overline{U}_{k}$ contains two boundary arcs $\alpha_1$, $\alpha_2\subset \partial D^{2}$ and $\partial U_{k}\setminus (\alpha_1\cup \alpha_2)=\beta_1\cup \beta_2$, where $\beta_i \subset \varphi(F)$, $i = 1, 2$.

The set $\Cl{U}_{k} \setminus (\beta_{1} \cup \beta_{2})$ divides a disk $D^{2}$ hence arcs $\beta_{1}$ and $\beta_{2}$ do not belong to the image of the same tree of $F$. Suppose that $\beta_{1} \subseteq \varphi(T_{i})$, $\beta_{2} \subseteq \varphi(T_{j})$, $i \neq j$.

It is obvious that any arc $\alpha_{1}$, $\alpha_{2}$ does not contain other images of vertices of $F$ besides its endpoints. By using $A3$ we can conclude that the only images of vertices of $G$ that are contained in $\alpha_{1}$ and $\alpha_{2}$ are their endpoints. Let us denote by $y_{rs}$, $r, s \in \{1, 2\}$, a common endpoint of $\alpha_{r}$ and $\beta_{s}$.

 From Lemma~\ref{main_lemma} it follows that $f(\beta_{1}) = g(T_{i}) = c_{i} \in \rr$, $f(\beta_{2}) = g(T_{j}) = c_{j} \in \rr$, and $f$ has no local extrema on arcs $\alpha_{r} \setminus \{y_{r1}, y_{r2}\}$, $r \in \{1, 2\}$. Therefore $c_{i} \neq c_{j}$ and maps $\psi_{1} : \alpha_{1} \rightarrow [0, 1]$, $\psi_{2} : \alpha_{2} \rightarrow [0, 1]$,
\begin{align*}
\psi_{1}(z) & = \frac{f(z)-f(y_{11})}{f(y_{12})-f(y_{11})} = \frac{f(z)-c_{i}}{c_{j}-c_{i}} \,, \\
\psi_{2}(z) & = \frac{f(z)-f(y_{21})}{f(y_{22})-f(y_{21})} = \frac{f(z)-c_{i}}{c_{j}-c_{i}} \,,
\end{align*}
are homeomorphisms, moreover $\psi_{1}(y_{11}) = \psi_{2}(y_{21}) = 0$, $\psi_{1}(y_{12}) = \psi_{2}(y_{22}) = 1 \,.$

\begin{figure}[htbp]
\begin{center}
\includegraphics{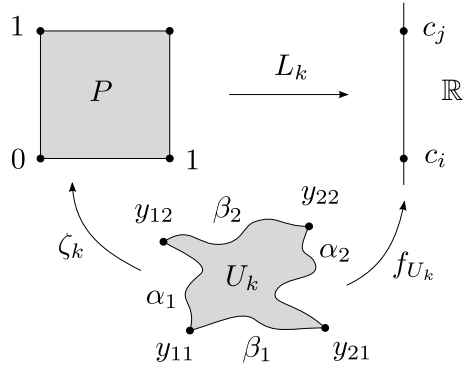}
\caption{Function on a connected domain with two boundary arcs.}\label{fig_case_2}
\end{center}
\end{figure}

We consider a set $P = [0, 1] \times [0, 1]$ and a map $L_{k} : P \rightarrow \rr$,
\[
L_{k} : (x, y) \mapsto c_{i} (1-y) + c_{j} y \,.
\]

Let $\eta_{s} : \beta_{s} \rightarrow [0, 1]$ be homeomorphisms such that $\eta_{s}(y_{1s}) = 0$, $\eta_{s}(y_{2s}) = 1$, $s \in \{1, 2\}$.

Let us consider a map $\zeta_{k}^{0} : \partial U_{k} \rightarrow \partial \hat{P}$,
\[
\zeta_{k}^{0}(z) =
\left\{
\begin{array}{ll}
	(0, \psi_{1}(z))\,, & \mbox{for } z \in \alpha_{1} \,, \\
	(1, \psi_{2}(z))\,, & \mbox{for } z \in \alpha_{2} \,, \\
	(\eta_{1}(z), 0)\,, & \mbox{for } z \in \beta_{1} \,, \\
	(\eta_{2}(z), 1)\,, & \mbox{for } z \in \beta_{2} \,.
\end{array}
\right.
\]
It is easy to see that $\zeta_{k}^{0}$ is homeomorphism. By using  Schoenflies's theorem~\cite{Newman} we can extend the homeomorphism $\zeta_{k}^{0}$ to a homeomorphism $\zeta_{k} : \Cl{U}_{k} \rightarrow \hat{P}$.

Let us consider a continuous function, see Fig.~\ref{fig_case_2}
\[
f_{U_{k}} = L_{k} \circ \zeta_{k} : \Cl{U}_{k} \rightarrow \rr \,.
\]
It is evident that this function locally  can be represented as a projection on coordinate axis in all points of $\Cl{U}_{k} \setminus \{y_{11}, y_{12}, y_{21}, y_{22}\}$.

By analogy with case 1, we prove that $\restrict{f_{U_{k}}}{\partial U_{k}} = \restrict{f}{\partial U_{k}}$.

\medskip

The union of $\{\Cl{U}_{k}\}$ generates a finite closed cover of $D^{2}$. In addition, it follows from the construction that if $z \in \Cl{U}_{k} \cap \Cl{U}_{s}$ for some $k \neq s$ then $z \in \varphi(G)$ and $f_{U_{k}}(z) = f_{U_{s}}(z) = f(z)$. Therefore we can extend a function $f$ from the set $\varphi(G)$ into $D^{2}$ with the help of the following relation
\[
f(z) = f_{U_{k}}(z) \,, \quad \mbox{for } z \in \Cl{U}_{k} \,.
\]
The cover $\{\Cl{U}_{k}\}$ is fundamental thus the function $f : D^{2} \rightarrow \rr$ is continuous.

$G$ will be indentified with its image $\varphi(G)\subseteq D^2$ in the following discussion.


Let $T_{k}$ be a tree of $F$. Let us denote by $\Theta_{k}$ an union of domains of $\Theta = D^{2} \setminus \varphi(F)$ which are adjoined to $T_{k}$.
\begin{align*}
\Theta_{k} & = \bigcup_{j=1}^{m(k)} U_{j}^{k} \,, \\
\{ U_{1}^{k}, \ldots, U_{m(k)}^{k} \} & = \{ U_{j_1}, \ldots, U_{j_{m(k)}} \} \,.
\end{align*}

It should be noted that for any domain $U_{j}$ and arc $\alpha = T_{k} \cap \partial U_{j}$ we have in the first place $f(\alpha) = const = c_{k}$, secondly, either $f(z) > c_{k}$ for any $z \in U_{j}$ or $f(z) < c_{k}$ for any $z \in U_{j}$. Thus every domain $U_{j}^{k}$, $j \in \{1, \ldots, m(k)\}$ of $\Theta_{k}$ can be associated with sign either ``$+$'' or ``$-$'' depending on the sign of difference $f(z) - c_{k}$, $z \in U_{j}^{k}$.

It is easy to see that arcs of $\partial D^{2}$ connecting the images of adjacent vertices of $G$ are connected components of the set $\Gamma_{k} = (\overline{\Theta}_{k} \cap \partial D^{2}) \setminus \varphi(V^{\ast})$. Therefore from Definition~\ref{ozn_cr_cycle} and Lemma~\ref{main_lemma} it follows that $f$ is monotone on any arc of $\Gamma_{k}$. By definition of $\Gamma_{k}$ and $A3$ exactly one of endpoints of any arc of $\Gamma_{k}$ is an image of vertex of tree $T_{k}$. Thus every arc $S$ of $\Gamma_{k}$ can be associated with a sign either ``$+$'' or ``$-$'' depending on the sign of difference $f(z) - c_{k}$, $z \in S$.

Let us prove that in a neighborhood of any vertex $\varphi(v)$, where $v\in V_k\subset T_k$, the signs of domains, whose boundaries are the images of edges adjacent to $v$ alternate. We should remark that for every vertex $v$ of $ V_{ter}(T_k)$ this follows from $A3$.

Let $\varphi(v)$ be a vertex such that $v\in V_k\setminus V_{ter}(T_{k})$. Suppose that in a neighborhood of some point of $e_m\backslash V(G)$ there exist two domains $U_m$ and $U_{m+1}$ which are adjoint to the edge $e_m$  (which is adjacent to $v$) such that they have the same sign. Then from $A2$ it follows that there exist an other edge $e_n$ that is adjacent to $v$ and both of its adjoining domains $U_n$ and $U_{n+1}$ have the same sign, cases $U_{m+1}=U_n$, $U_m=U_{n+1}$ are not excluded. Since for an edge $e_n$ ($e_m$) there exist one more vertex $v'_1$ ($v''_1$) which is adjacent to it then by analogy for the vertex $v'_1$ ($v''_1$) we can find a vertex $v'_2$ ($v''_2$) such that it has two adjacent edges adjoining to domains with the same sign and so on. Tree is finite, so for sequence of vertices $v,v'_1,v'_2,\ldots,v'_{s_1}$ ($v,v''_1,v''_2,\ldots,v''_{s_2}$) there exists a vertex $v'_{s_1}$ ($v''_{s_2}$) such that $v'_{s_1}\in V^{\ast}_k$ ($v''_{s_2}\in V^{\ast}_k$).

We consider two types of vertices from $V^{\ast}_k\setminus V_{ter}(T_k)$.

\medskip
(i) if $deg(v'_{s_1})=2k+1>3$, then a number of domains adjoining to the edges which are incident to it is even.

Let us consider the following binary relation $\rho$ on the set of such domains. We will say that $V' \rho V''$ if domains $V'$ and $V''$ adjoin to a common edge $e$ which is incident to vertex $v'_{s_1}$ and going around $v'_{s_1}$ across the edge $e$ in positive direction we pass from $V'$ to $V''$.

It is easy to show that the relation $\rho$ is convenient and all domains generate $\rho$-chain, where its first and last elements are domains whose boundary contain arcs $S'$, $S'' \subseteq \Gamma_{k}$ adjoining to $v'_{s_1}$.

From $A3$ it follows that arcs $S'$ and $S''$ have different signs therefore first and last element of $\rho$-chain have the different signs.

From this and the fact that $\rho$-chain has even number of elements it follows that a number of its pairs of adjacent elements which have the same sign is even. We can apply our previous argument and add one more vertex $v_{s_1+1}'$ to the sequence of vertices $v$, $v_{1}', \ldots, v_{s_1}'$.

\medskip
(ii) if $deg(v'_{s_1})=2k>3$, then a number of domains adjoining to edges which are adjacent to it is odd.

Just as in (i) we consider the relation $\rho$ on the set of such domains and order them into $\rho$-chain.

Contrary to the previous case a length of $\rho$-chain is odd and by Condition $A3$ first and last its elements have the same sign. Similarly, in this case a number of pairs of adjacent elements of $\rho$-chain which have the same signs is even. Therefore we can add one more vertex  $v_{s_1+1}'$ to the sequence of vertices $v$, $v_{1}', \ldots, v_{s_1}'$.

\medskip

From the finiteness of tree it follows that there exists a vertex $v'$ ($v''$) such that $v'\in V_{ter}(T_k)$ ($v''\in V_{ter}(T_k)$) and an edge which is incident to it adjoins to domains with the same signs. Thus arcs $S'_l$ and $S'_p$ with the endpoint $v'$ ($v''$) have the same sign but it contradicts to $A3$.

Let us consider the restriction of $f$ to $\partial D^{2}$. Local extrema of $f$ are points $\varphi(v_{i})$ corresponding to vertices $v_{i}$ such that $v_{i}\in \gamma$ and $deg(v_{i})=2k$, see Lemma~\ref{main_lemma}. From the finiteness of $G$ follows the finiteness of number of local extrema on $\partial D^{2}$.


\end{proof}

Let $G$ be $\Delta$--graph. From Theorem~\ref{main_theorem} it follows that there is a pseudoharmonic function $f$ on disk which corresponds to a graph $G$. But, in general this function is not uniquely defined since we in no way restrict the choice of a monotone map $g : G \rightarrow \rr$. Thus for non comparable vertices $v'$ and $v''$ of graph $G$ the relation $g(v')= g(v'')$ is not necessarily valid.

It is easy to construct an example of $\Delta$--graph $G$ and two monotone maps $g_{1}, g_{2} : G \rightarrow \rr$ which satisfy Lemma~\ref{main_lemma} but for some pair of non comparable vertices $v'$, $v'' \in V(G)$ the following correlations hold true $g_{1}(v') < g_{1}(v'')$ and $g_{2}(v') > g_{2}(v'')$.

\begin{theorem}
Let a graph $G$ be $\Delta$--graph.

$G$ satisfies Condition $A4$ iff a strict partial order of a graph $G$ coincides with a strict partial order of a diagram $P(f)$ of some pseudoharmonic function $f$ that corresponds to $G$.
\end{theorem}

\begin{proof}
Let $P(f)$ be a combinatorial diagram of some pseudoharmonic function $f$. We remind that a partial order on vertices of $P(f)$ is induced by a function $f$ with the help of the following relation
\[
v' < v'' \,, \quad \mbox{if } f \circ \psi(v') < f \circ \psi(v'') \,, \quad v', v'' \in V(P(f)) \,.
\]
We note that vertices $v'$ and $v''$ are non comparable iff their images are on the same level set of $f$. Hence a graph $P(f)$ satisfies $A4$.

\medskip

Suppose that $G$ satisfies $A4$. The binary relation ``to be non comparable'' on the set of vertices $V(G)$ of $G$ is transitive, symmetric and reflexive. So, in the proof of Lemma~\ref{main_lemma} we can consider instead of $\ff$ a partition $\tilde{\ff}$  whose elements are classes of non comparable elements with regards to the order on $V(G)$. Then due to condition A4 the projection $\tilde{\pi} : V(G) \rightarrow V(G)/\tilde{\ff}$ induces a relation of partial order on quotient space $\tilde{V} = V(G)/\tilde{\ff}$ such that every pair of elements $\tilde{v}'$, $\tilde{v}'' \in \tilde{V}$ is comparable. Therefore the partially ordered space $\tilde{V}$ is linearly ordered and every monotone map $\tilde{g} : \tilde{V} \rightarrow \rr$ is isomorphism onto its image. A map $g = \tilde{g} \circ \tilde{\pi}$ satisfies the condition that any pair of vertices $v'$, $v'' \in V(G)$ is non comparable iff $g(v') = g(v'')$.

In the same way as in Lemma~\ref{main_lemma} we extend the function $g$ on $G$ and use this extension to construct a pseudoharmonic function $f$. By the construction the partial order induced on $V(G)$ by $f$ is the same as the original partial order on $V(G)$.
\end{proof}



\begin{thebibliography}{99}

\bibitem{Kp}
\textit{Kaplan}
Topology of level curves of harmonic functions.
Transactions of Amer. Math. Society -- vol. 63, N 3 (1948) -- pp.~514--522.

\bibitem{Yu}
\textit{Yurchuk I.}
Topological equivalence of functions of class $F(D^{2})$. (in Ukrainian)
Zb. prac Inst. Math. NAS Ukraine, 2006. -- V.3, N 3. --  474--486 pp.

\bibitem{Bo}
\textit{Boothby W.M.}
The topology of regular curve families with multiple saddle points.
Amer. J. Math. -- 1951. -- vol. \textbf{73}. -- pp.~405--438.

\bibitem{Mrt}
\textit{Morse M.}
The topology of pseudo-harmonic functions.
Duke Math. J. -- 1946. -- vol. \textbf{13}. -- pp.~21--42.

\bibitem{PY1}
\textit{Polulyakh E., Yurchuk I.}
On the conditions of topological equivalence of pseudoharmonic functions defined on disk.
arXiv:0910.3647v1 [math.GN]

\bibitem{PY}
\textit{Polulyakh E., Yurchuk I.}
On the criteria of D-planarity of a tree.
arXiv:0904.1367v1 [math.GN]

\bibitem{Nov}
\textit{Novak V.}
Cyclically ordered sets.
Czechoslovak Math. Journal., 1982. -- vol. 32(107). -- pp.~460--473.

\bibitem{Mel}
\textit{Mel'nikov O. V., Remeslennikov, V. N., Roman'kov, V. A., Skornyakov, L. A., Shestakov, I. P.}
General algebra. Vol. 1 (Russian)
Mathematical Reference Library, ``Nauka'', Moscow, 1990. 592~pp.

\bibitem{Newman}
\textit{Newman M. H. A.}
Elements of the topology of plane sets of points.
Cambridge: Cambridge Univ. Press, 1964 -- 214~pp.

\bibitem{RF}
\textit{Fuks D. B., Rokhlin V. A.}
Beginner's course in topology. Geometric chapters.
Translated from the Russian by A. Iacob.
Universitext. Springer Series in Soviet Mathematics. Springer-Verlag, Berlin, 1984. xi+519 pp.

\end{thebibliography}
\end{document}